\documentclass [11pt]{article}
\usepackage{amsmath, amsthm, amssymb}
\usepackage{graphicx}
\newfont{\bbb} {msbm10}
\newcommand{\R}{\Bbb{R}}
\newcommand{\bS}{\Bbb{S}}
\newcommand{\bB}{\Bbb{B}}
\newcommand{\T}{\Bbb{T}}

\newcommand{\qH}{\Bbb{H}}

\newcommand{\sbs}{\subset}
\newcommand{\ra}{\rightarrow}

\newcommand{\bphi}{{\bar{\phi}}}
\newcommand{\br}{{\bar{r}}}
\newcommand{\bt}{{\bar{t}}}

\newcommand{\p}{\partial}
\newcommand{\barp}{\bar{\partial}}
\newcommand{\pt}{\frac{\partial}{\partial t}}

\newcommand{\cE}{{\cal{E}}}

\newcommand{\B}{\Bbb{B}}

\newcommand{\HH}{\Bbb{H}}

\newcommand{\0}[1]{_{_{#1}}}

\usepackage[all]{xy}

\begin{document}

\title{Hyperbolic Extensions and Metrics $\epsilon$-Close to Hyperbolic}
\author{Pedro Ontaneda\thanks{The author was
partially supported by a NSF grant.}}
\date{}

\maketitle

\begin{abstract} We define the {\it Hyperbolic Extension} of
a Riemannian manifold {\it with a center}, and give
some properties of it. Our main result says to what
extent the hyperbolic extension of $M$ is {\it close to being
hyperbolic}, if we assume $M$ to be close to hyperbolic.

The results in this paper are used in the problem of smoothing
Charney-Davis strict hyperbolizations \cite{ChD}, \cite{O}.
\end{abstract}

\noindent {\bf \large  Section 1. Introduction.}

 Recall that hyperbolic $n$-space $\HH^n$ is isometric to $\HH^{k}\times \HH^{n-k}$ with
warped metric $(\cosh^2\, r)\,\sigma\0{\HH^{k}}+\sigma\0{\HH^{n-k}}$, where $\sigma\0{\HH^{l}}$ denotes the hyperbolic metric of
$\HH^{l}$, and  $r:\HH^{n-k}\ra[0,\infty)$ is the distance to a fixed point in $\HH^{n-k}$.
For instance, in the case $n=2$, since $\HH^1=\R^1$ we have that $\HH^2$ is isometric to $\R^2=\{(u,v)\}$ with warped metric $\cosh^2v\, du^2+ dv^2$. In the following paragraph we give a generalization of this
construction.\vspace{.1in}

Let $(M^n,h)$ be a complete Riemannian manifold with {\it center} $o=o_{_M}\in M$, that is, the exponential map $exp_o:T_oM\ra M$ is a diffeomorphism.
The warped metric 
\begin{center}$g=(\cosh^2 r)\, \sigma\0{\HH^k}+h$
\end{center} on $\HH^k\times M$
is the {\it hyperbolic extension (of dimension $k$)} 
of the metric $h$. Here $r$
is the distance-to-$o$ function on $M$.
We write $\cE_k(M)=(\HH^k\times M,g)$, and $g=\cE_k(h)$.
We also say that $\cE_k(M)$ is the {\it hyperbolic extension
(of dimension $k$) of $(M,h)$} (or just of $M$).
Hence, for instance, we have $\cE_k(\HH^l)=\HH^{k+l}$.
For $S\sbs M$ we write $\cE_k(S)=\HH^k\times S\sbs\cE_k(M)$.\vspace{.1in}

Also write $\HH^k=\HH^k\times \{o\0{M}\}\sbs\cE_k(M)$ and
we have that any $p\in\HH^k$ is a center of $\cE_k(M)$
(see Section 1). Since $\cE_k(\HH^l)=\HH^{k+l}$ one would
expect that if $M$ is, in some sense, close to $\HH^l$, then
$\cE_k(M)$ would be close to $\HH^{k+l}$.
As mentioned in the abstract our main result states to what
extent the hyperbolic extension of $M$ is close to being
hyperbolic, if we assume $M$ to be close to hyperbolic.
Our definition of ``close to hyperbolic'' is a chart-by-chart
definition and it is given in the next paragraph (see Section 2 for more details). In this paper we also give some properties
of hyperbolic extensions and introduce a set of coordinates
well suited to study these objects. 
\vspace{.1in}

Let $\bB\sbs\R^{l-1}$ be the unit $(l-1)$-ball, 
with the flat metric $\sigma\0{\R^{l-1}}$.  Write  $I_\xi=(-1-\xi,1+\xi)$, $\xi>0$.
Our basic models are $\T_\xi=\bB\times I_\xi$,  with hyperbolic metric $\sigma=e^{2t}\sigma\0{\R^{l-1}}+dt^2$. 
The number $\xi$ is called the {\it excess} of $\T_\xi$.
(The reason for introducing $\xi$ will become clear in 
Theorem A below, see Remark 1).
Let $(N^l,g)$ be a Riemannian manifold and $S\sbs N$. We say that $S$ is $\epsilon$-{\it close to hyperbolic} if there is $\xi>0$ such that for every $p\in S$ there is an {\it $\epsilon$-close to hyperbolic 
chart with center $p$}, that is, there is a chart
$\phi :\T_\xi\ra N$, $\phi(0,0)=p$,  such that $|\phi^*g-\sigma|\0{C^2}<\epsilon$.  
Here $|.|_{C^2}$ is the $C^2$ norm (see Section 2)The number $\xi$ is called the {\it excess } of the charts
(which  is fixed). \vspace{.1in}

Let $(N^l,g)$ have center $o$.
Using the exponential map $exp_o$ we shall sometimes identify $N$ with
$\R^l$, and  $N-\{ o\}$ with $\bS^{l-1}\times\R^+$. Let $S\sbs N $.
We shall say that $g$ is {\it radially $\epsilon$-close to hyperbolic on
$S$ (with respect to $o$)} if, 
for every $p\in S$ there is an  $\epsilon$-close to hyperbolic 
chart $\phi$ with center $p$ and, in addition, the chart 
$\phi$ respects the product structure of $\T_\xi$
and $N-o=\bS^n\times\R^+$, that is $\phi(. , t)=(\phi_1(.), t+a)$, for some constant $a$ depending only on the
$\phi$ (see Section 2  for details). Here the ``radial" directions
are $(-1-\xi,1+\xi)$ and $\R^+$ in $\T_\xi$ and $N-o$, respectively.
Sometimes we will just say ``$S$ is radially $\epsilon$-close to hyperbolic" when it is clear from the context which metric is 
being considered.
\vspace{.1in}

\noindent{\bf Remark 1.1} Our definition of
radially $\epsilon$-close to hyperbolic metric is well suited
to study metrics away from the center, but near the center 
this definition is not useful. This is because: (1) the need for some space to fit the charts, and (2) the form of our specific fixed model $\T_\xi$.
For instance hyperbolic $n$-space
$\HH^{n}$ is {\sf not} radially $\epsilon$-close to hyperbolic
near a (chosen) center. In fact 
$\HH^n$ is radially $\epsilon$-close to hyperbolic outside
the ball $\B_{\sf a}(\HH^n)$ of radius ${\sf a}$, for a constant  ${\sf a}={\sf a}(\epsilon,n,\xi)$ depending only on
$\epsilon$, $n$, $\xi$ (see 4.14 in \cite{O1}).\vspace{.1in}

The next Theorem states that if $S\sbs M$ is radially $\epsilon$-close
to hyperbolic then $\cE_k(S)$ is $\eta$-close to hyperbolic,
where $\eta$ depends on $\epsilon$ and the distance $r$ to $\HH^k\sbs\cE_k(M)$.  We use the notation $B_a=\B_a(M)\sbs M$ for the ball of radius $a$ centered at the center $o=o\0{M}$. We assume $\xi>0$.\vspace{.1in}

\noindent {\bf Theorem A.} {\it Let $M^n$ have center $o$, and
$S\sbs M-\{ o\}$. Assume  $S$ is radially $\epsilon$-close to hyperbolic, with charts of excess $\xi$. Then $\cE_k(S- B_a)$
is radially $\eta$-close to hyperbolic 
(with respect to any point in $\HH^k$), provided}
$$
C\,\Big(\,\epsilon\,+\, e^{-a}\,\Big)\,\,\leq\,\, \eta
$$
\noindent {\it where $C$ is a constant depending on $k$ and $\xi$.  Moreover,  $\cE_k(S- B_a)$
is radially $\eta$-close to hyperbolic, with charts of excess $\xi'$, provided $$0\,<\,\xi'\,<\,\xi\,-\,L e^{-a}$$
\noindent where $L$ is a constant depending on $k$ and $\xi$.}
\vspace{.2in}

\noindent {\bf Remarks.} 

\noindent {\bf 1.} Note that the excess of the charts decreases.
This is the main reason to introduce the excess. In \cite{O2} we described another
geometric process, warp forcing, which also reduces the excess of the charts.

\noindent {\bf 2.} In the statement of the theorem above 
$S$ is radially $\epsilon$-close to hyperbolic with respect to the decomposition $M-\{ o\}=\bS^{n-1}\times\R^+$ and
$g$ is radially $\eta$-close to hyperbolic  with respect to the decomposition $\cE_k(M)-\{ o\}=\bS^{n+k-1}\times\R^+$.

\noindent {\bf 3.} An explicit formula for $L$ is given at the end
of Appendix A. It is implicit in Theorem A
that $r$ cannot be too small, that is, we want $\xi\,-\,L e^{-r}\geq 0$, hence we want $r\geq ln(L)-ln(\xi)$.

\noindent {\bf 4.} We can take $C=C(k,\xi)=
2(2+3\xi+\xi^2)e^{^{1+\xi}}L\, +\, C_4$, where $C_4=C_4(k,\xi)=C_1(c\0{\bS^k},k,\xi)$, 
with $C_1$ as in Corollary 4.3 of
\cite{O1}, and $c\0{\bS^k}$ is such that $\bS^k$ is 
$c\0{\bS^k}$-bounded (see Section 4 of \cite{O1}).
\vspace{.1in}

We now deal with a natural and useful class of metrics.
These are metrics on $\R^n$ (or on a manifold with center)
that are already hyperbolic on the ball
$B_{a}=B_a(0)$ of radius $a$ centered at $0$, and are radially $\epsilon$-close to hyperbolic outside
$B_{a'}$
(here $a'$ is slightly less than $a$). Here is the detailed
definition. Let $M^n$ have center $o$ and
let $B_{a}=B_a(o)$ be the ball on $M$ of radius $a$ centered at $o$. We say
that a metric $h$ on $M$ is $(B_a,\epsilon)$-{\it close to hyperbolic, with charts of excess $\xi$}, \, if \vspace{.1in}

\begin{enumerate}
\item[  (1)]  On $B_{a}-\{o\}=\bS^{n-1}\times (0,a)$
we have $h=\sinh^2(t)\sigma\0{\bS^{n-1}}+dt^2$. Hence $h$
is hyperbolic on $B_a$.
\item[(2)]  the metric $h$ is 
radially $\epsilon$-close to hyperbolic outside $B_{a-1-\xi}$,
with charts of excess $\xi$.
\end{enumerate}\vspace{.1in}

\noindent {\bf Remarks.}

\noindent {\bf 1.} We have dropped the word ``radially" 
to simplify the notation. But it does appear in condition (2),
where ``radially" refers to the center of $B_a$.

\noindent{\bf 2.} We will always assume $a>{\sf a}+1$,
where {\sf a} is as in 1.1. Therefore conditions (1), (2) 
and Remark 1.1 imply
a stronger version of (2):\vspace{.1in}

  (2') the metric $h$ is 
radially $\epsilon$-close to hyperbolic outside $B_{\sf a}$,
with charts of excess $\xi$.\vspace{.1in}

\noindent This is the reason why we demanded radius $a-1-\xi$
in (2), instead of just $a$. 
 \vspace{.1in}

Metrics that are $(B_a,\epsilon)$-close to
hyperbolic are very useful, and are key objects in
\cite{O}. See also \cite{O1}, \cite{O4}.
Our next result answers the following question:
\vspace{.1in}

\noindent {\bf Question.} {\it What can we say about the hyperbolic extension of a $(B_a,\epsilon)$-close to
hyperbolic metric?}\vspace{.1in}

\noindent {\bf Theorem B.}  {\it  Let $M^n$ have center $o$. Assume $M$ is $(B_a,\epsilon)$-close to hyperbolic, with charts of excess $\xi>0$. Then $\cE_k(M)$ is $(B_a,C_2\epsilon)$-close to hyperbolic, with charts of excess $\xi'$,
provided $a$ is sufficiently large. Explicitly we want
\newline\hspace*{2.5in}$a\,\geq \,R\,=\,R(\epsilon,k,\xi)$
\vspace{.05in}

\noindent Here $C_2=C_2(n,k,\xi)$, and  \,$\xi'=\xi-e^{-a/2}>0$.}\vspace{.1in}

\noindent {\bf Remark.} 
The constant  $R$ is defined as
$R=ln\big(\frac{1}{\epsilon}\big)+ln(L)+1+\xi$.
Here $C_2=C_2(n,k,\xi)=C_1'e^{1+\xi}+C$, where $C_1'$ is as in 4.13 of \cite{O1}, and $C$, $L$ are as in Theorem A.\vspace{.2in}

The results in this paper are used to construct negatively curved Riemannian smoothings of Charney-Davis strict
hyperbolizations of manifolds \cite{ChD}, \cite{O}. In the next paragraph we give an idea how the objects and results in this paper is used in \cite{O}.
\vspace{.1in}

In the same way that a cubical complex is made of basic pieces (the cubes
$\square^k$),
the hyperbolization $h(K)$ of a cubical complex $K$ is 
also made of
basic pieces: pre-fixed hyperbolization pieces $X^k$. Indeed one begins with a cubical complex $K$ and replaces each cube of dimension $k$
by the hyperbolization piece of the same dimension. Cube complexes have
a piecewise flat metric induced from the flat geometry of the cubes.
Likewise the Charney-Davis hyperbolizations have a piecewise hyperbolic
structure because the Charney-Davis hyperbolization pieces
are hyperbolic manifolds (compact, with boundary and corners).
To see how singularities appear one can first think about the manifold 2-dimensional 
cube case. If $K^2$ is a 2-dimensional manifold cube complex then
its piecewise flat metric is Riemannian outside the vertices. A vertex is
a singularity if and only if the vertex does not meet exactly four cubes.
The picture is exactly the same for $h(K^2)$. 
These point singularities in $h(K^2)$ can be smoothed out
easily using warping methods.
In higher dimensions
the singularities of $K^n$ and $h(K)$ appear in (possibly the whole of)
the codimension 2 skeletons $K^{(n-2)}$ and $h(K^{(n-2)})$, respectively.
In \cite{O} the idea of smoothing the piecewise hyperbolic metric on
$h(K)$ is to do it inductively down the dimension of the skeleta.
One begins with the $(n-2)$-dimensional pieces $X^{n-2}$. Transversally
to each $X^{n-2}$ (that is, on the union of geodesic segments emanating
perpendicularly to $X^{n-2}$, from a fixed point in $X^{n-2}$) one has
essentially the 2-dimensional picture mentioned above. 
Once we
solve this transversal problem we extend this transversal smoothing by taking a warp product
with $X^{n-2}$; this product method is the {\it hyperbolic extension}
process studied in this paper. After applying this method
we have a smoothing on a (tubular) neighborhood of the piece $X^{n-2}$.
Caveat: we do not want to actually have a smoothing on a neighborhood
of the whole of $X^{n-2}$, since we will certainly have 
matching problems for different $X^{n-2}$ meeting on a common
$X^{n-3}$; so we only want a smoothing on a neighborhood of the $Z^{n-2}$,
where $Z^{n-2}\sbs X^{n-2}$ is just a bit ``smaller" than $X^{n-2}$,
so that the neighborhoods of the $Z^{n-2}$ are all disjoint. Next 
step is to smooth around the $X^{n-3}$ (or, specifically the $Z^{n-3}$).
The metric is already smooth outside a neighborhood of the $(n-3)$-skeleton. Transversally to each $X^{n-3}$ we
have a 3 dimensional problem.
(It helps to have a 3 dimensional picture in mind, like in dimension 2).
It happens that if we did things with care in the first step (around the
$Z^{n-2}$) the metric in the 3 dimensional transversal problem
is radially $\epsilon$-close to hyperbolic outside some large ball B.
At this point we use the method of {\it hyperbolic forcing} introduced in
\cite{O4} to
 to extend the metric to a Riemannian metric
on the ball B, getting rid, in this way, of the transverse
singularity. The resulting
metric is still radially $\eta$-close to hyperbolic, with an $\eta$
that can be controlled. Once the transversal 3 dimensional problem
is solved we extend this smoothing to neighborhoods of the $Z^{n-3}$
using hyperbolic extension. Next we do the same for the $Z^{n-4}$
and so on.  Along the whole process we want to control
how far our metrics are from being hyperbolic; this is why
the main theorems in this paper are fundamental in \cite{O}.\vspace{.2in}

Here is a brief description of the paper. In Section 2 we give more details about the 
$C^2$ norm.
In Section 3 we give some basic facts about hyperbolic extensions and give some
properties. In Section 4 we introduce some useful coordinates
in hyperbolic extensions. In Section 5 we prove
Theorems A and B. There are 2 appendices in which we deal with some technical details.\vspace{.3in}

\noindent {\bf \large  Section 2. The $C^2$-Norm.}\\ Let $A\sbs \R^n$ be an open set.
Let $|.|\0{C^2(A)}$ denote the uniform $C^2$-norm of $\R^l$-valued functions on $A$, i.e. if $f=(f\0{1},...,f\0{l}):A\ra \R^l$,
then $|f|\0{C^2(A)}=sup\0{z\in A,
\,\,1\leq i\leq l,\,\,1\leq j,k\leq n}\{ |f\0{i}(z)|, |\p\0{j}f\0{i}(z)|, |\p\0{j,k}f\0{i}(z)|\}$. Sometimes we will write $|.|\0{C^2}=|.|\0{C^2(A)}$ when the context is clear.
Given a Riemannian metric $g$ on $A$, the number $|g|\0{C^2(A)}$ is computed considering $g$ as the $\R^{n^2}$-valued function $z\mapsto (g_{ij}(z))$ where, as usual,
$g\0{ij}=g(e_i,e_j)$, and the $e_i$'s are the canonical vectors in $\R^{n}$. 
\vspace{.1in}

The $C^2$-norm $|.|\0{C^2}$ mentioned in the  
definition of an  $\epsilon$-close to hyperbolic Riemannian manifold in the Introduction is $|.|\0{C^2}=|.|\0{C^2(\T_\xi)}$. 
If $(M,g)$ is $\epsilon$-close to hyperbolic (or radially $\epsilon$-close 
to hyperbolic) we will also say that the metric $g$ is $\epsilon$-close to hyperbolic (or radially $\epsilon$-close to hyperbolic).
\vspace{.3in}

\newpage

\noindent{\bf \large Section 3. Hyperbolic Extensions.}

Let $M^n$ be a complete Riemannian manifold with {\it center} $o=o_{_M}\in M$, that is, the exponential map $exp_o:T_oM\ra M$ is a diffeomorphism.
In particular $M$ is diffeomorphic to $\R^n$. 
For instance if $M$ is Hadamard manifold every point is a center point.
Denote the metric on $M$ by $h$.\vspace{.1in}

 In this paper we will use the same symbol ``$o$" to denote a center of a Riemannian manifold unless it is necessary to
specify the manifold, in which case we will write $o_{_M}$ if $o$ is a center of $M$.\vspace{.1in}

Let $r:M\ra [0,\infty)$ be the distance to $o$. Then 
\begin{enumerate}
\item[{\bf i.}] we have that $r(exp_o v)=h_o(v,v)^{1/2}$, hence $r$ is continuous and smooth on $M-\{ o\}$. Also $r^2$ is smooth on $M$.

\item[{\bf ii.}]  The (images of the) geodesic rays $exp_o(\R^+ v)$ are convex sets in $M$, and the geodesics lines $exp_o(\R v)$ are totally geodesic in $M$.
Here $\R^+= (0,\infty)$.

\item[{\bf iii.}]  the function $dr$ is strictly distance decreasing on non-radial vectors. That is, for $v\in TM-T_oM$ we have $|dr(v)|\leq h(v,v)^{1/2}$ and $|dr(v)|= h(v,v)^{1/2}$ if and only if $v$ is radial, i.e. tangent to a geodesic passing through $o$. (This follows from the Gauss Lemma and the fact that $r\circ (exp_o)^{-1}:T_oM\ra\R$ is the
euclidean distance to the origin.)
\end{enumerate}

Using the diffeomorphism $exp_o$ onto $M$ and an identification
of $T_oM$ with
$\R^n$ via some fixed choice of an orthonormal basis in $T_oM$,
we can identify $M$ with $\R^n$ and $M-\{o\}$ with $\bS^{n-1}\times\R^+$.
Therefore
 the metric $h|_{M-\{ o\}}$ can be written 
as  $h_r+dr^2$ on $\bS^{n-1}\times \R^+$.
Also we shall call the set of ($h$-geodesic) rays
$t\mapsto (x,t)\in \bS^{n-1}\times \R^+$ the {\it ray structure of
$h$ with respect to o}.\vspace{.1in}

As mentioned in the Introduction the warped metric $g=(\cosh^2 r)\, \sigma\0{\HH^k}+h$ on $\HH^k\times M$
is a {\it hyperbolic extension} 
of the metric $h$ on $M$, and sometimes we will also write
$g=\cE_k(h)$.
Note that, even though $r$ is not smooth at $o$, the warping function $\cosh\,r$ is smooth on
$M$ because $\cosh$ is a smooth  even function.
Since $M$ is complete we have that $\cE_k(M)$ is also complete (see \cite{BisOn}, p.23).\vspace{.1in}

For instance, if  $M=\HH^{l}$ then the hyperbolic extension $\cE_k(\HH^{l})$ is 
hyperbolic $(k+l)$-space $\HH^{k+l}=\HH^k\times\HH^{l}$, with metric $(\cosh^2\, r)\,\sigma\0{\HH^{k}}+\sigma\0{\HH^{l}}$.\vspace{.1in}

For a subset $A\sbs\HH^k$ we shall write $\cE_A(M)= A\times M\sbs\cE_k(M)$, with
the metric $\cE_k(h)$ restricted to the set $\HH^k\times A$.\vspace{.1in}

We will write $\HH^k=\HH^k\times \{ o\}\sbs\cE_k(M)$. The hyperbolic extension, away from 
$\HH^k\sbs \cE_k(M)$ can be described in an alternative way: it is isometric to
$(\HH^k\times \bS^{n-1})\times (0,\infty)$ with metric $(\cosh^2r)\,\sigma\0{\HH^k}+h_r+dr^2$.\vspace{.1in}

It is known that every $\{y\}\times M$ is totally geodesic in 
$\cE_k(M)$ (see \cite{BisOn}, p.23).
Let $\eta$ be a complete geodesic line in $M$ passing though $o$
and let $\eta^+$ be one of its two geodesic rays (beginning at $o$) . Then $\eta$ is 
a totally geodesic subspace of $M$ and $\eta^+$ is convex (see item (ii) above). Also, let $\gamma$ be a complete geodesic line in $\HH^k$.
\vspace{.1in}

\noindent {\bf Lemma 3.1.} {\it 
We have that\,\, $\gamma\times \eta^+$ is a convex subspace of $\cE_k(M)$
and $\gamma\times \eta$ is totally geodesic in $\cE_k(M)$.}
\vspace{.1in}

\noindent {\bf Proof.} Let $\pi_\gamma:\HH^k\ra\gamma\sbs\HH^k$ denote the orthogonal projection, and note that
$d\pi_\gamma$ is distance non-increasing, i.e. $\sigma\0{\HH^k}(v,v)\geq \sigma\0{\HH^k}
(d\pi_\gamma(v),d\pi_\gamma(v))$, for $v\in T\HH^k$. Moreover, the equality holds if and only if $v\in T\gamma$.\vspace{.1in}

We assume $\eta:[0,\infty)\ra\eta\sbs M$ to be parametrized by the arc-length, that is, it is a speed-one geodesic ray.
Let $\pi_\eta:M\ra\eta^+$ denote the proper map $\pi_\eta(p)=\eta(r(p))$ . Note that $\pi_\eta$ is
smooth on $M-\{ o\}$ and  item (iii) above implies that $\pi_\eta$ is
strictly distance decreasing on non-radial tangent vectors on $M-\{ o\}$.  \vspace{.1in}

We prove that $\gamma\times \eta^+$ is convex.
 Let $\alpha:[0,1]\ra\HH^k\times M$ with end points on $\gamma\times \eta^+$.
Write 
$\alpha(u)=(a(u),b(u))\in\HH^k\times M$.  To prove that $\gamma\times \eta^+$ is convex
it is enough to prove that there is a curve $\beta$
on $\gamma\times \eta^+$,
with the same
end points as $\alpha$, but with length less or equal the length of
$\alpha$, and strictly less length if $\alpha$ is not contained in
$\gamma\times \eta^+$. For this, assume  for now that $b(u)\neq o$
for all $u\in [0,1]$, and
let $\beta=(\pi_\gamma\, a,\pi_\eta\, b)$. Hence
the length of $\alpha'=(a',b')$ is less or equal the length of $\beta'=
(d\pi_\gamma (a'),d\pi_\eta (b'))$. Therefore
the length of $\alpha$ is greater than the length of $\beta$, unless $a=\pi_\gamma (a)$ and $b$ is contained in a ray. And, by continuity the same holds
without the assumption that $b(u)\neq o$. Therefore $\gamma\times\eta^+$ is convex because $\beta$
is a path in $\gamma\times\eta^+$. \vspace{.1in}

We prove that $\gamma\times \eta$ is totally geodesic. Let $\eta^-=\overline{\eta-\eta^+}$ be the
``other" geodesic ray of $\eta$. Then $\gamma\times\eta^-$ is also convex.
Therefore $\gamma\times \eta-\gamma$ 
is totally geodesic hence the second fundamental 
of $\gamma\times\eta$ vanishes there. 
By continuity this form vanishes on the whole of $\gamma\times\eta$. 
This proves the lemma. \hfill $\Box$
\vspace{.1in}

Using the same construction in the proof of Lemma 3.1
one can show that every $\{y\}\times M$ is convex in 
$\cE_k(M)$. Moreover, since $r$ has a strict minimum at $o$,
one can also show that $\HH^k$ is convex in 
$\cE_k(M)$.\vspace{.1in}

\noindent {\bf Corollary 3.2.} {\it We have that \,$\HH^k\times\eta^+$  and $\gamma\times M$ are convex in $\cE_k(M)$. Also  \,$\HH^k\times\eta$
is totally geodesic in $\cE_k(M)$.}\vspace{.1in}

\noindent {\bf Proof.} For $\HH^k\times\eta$ just replace $\gamma$ by $\HH^k$ and $\pi_\gamma$ by the identity in the proof of
Lemma 3.1. For $\gamma\times M$ replace $\beta$ in the proof of Lemma 3.1 by  $\beta=(\pi_\gamma\, a,  b)$. 
This proves the Corollary.\hfill $\Box$\vspace{.1in}

\noindent {\bf Remarks 3.3.}  \\
\noindent {\bf 1.}  Note that $\HH^k\times\eta$ (with metric
induced by $\cE_k(M)$)
is isometric to $\HH^k\times \R$ with warped metric $\cosh^2 v\, \sigma\0{\HH^k}+dv^2$, which is just hyperbolic $(k+1)$-space $\HH^{k+1}$. Also $\gamma\times\eta$
is isometric to $\R\times \R$ with warped metric $\cosh^2 v\, du^2+dv^2$, which is just hyperbolic 2-space $\HH^2$.
In particular every point in $\HH^k=\HH^k\times\{ o\} \sbs\cE_k(M)$ is a center point.\\
\noindent {\bf 2.} It follows from Lemma 3.1 and Remark 1 that
the ray structure of $\cE_k(h)$ with respect to any center $o\0{\HH^k}\in\HH^k
\sbs\cE_k(M)$ only depends on the ray structure of $M$ and the center $o\0{\HH^k}$.\\
\noindent {\bf 3.} Denote by $\B_r(M)$ the ball on $M$ of radius $r$ 
about the center.
Note that if $h$ and $h'$ on $M$ have the same ray structures
then the balls $\B_r(M)$ coincide.\\
\noindent {\bf 4.} Recall that $\HH^k$ is convex in $\cE_k(M)$.
Moreover, for $l\leq k$, we also have $\HH^l\sbs\HH^k\sbs\cE_k(M)$ is
convex. If $h$ and $h'$ on $M$ have the
same ray structures then the $r$-neighborhoods (with respect to
$h$ and $h'$) of the convex subset
$\HH^l$ coincide.

\vspace{.3in}

\noindent {\bf \large Section 4. Coordinates on $\cE_k(M)$.} \\
Recall that we are identifying $M-\{ o\}$, $o=o\0{M}$, with $\bS^{n-1}\times \R^+$, and sometimes we shall denote a point
$v=(u,r)\in \bS^{n-1}\times\R^+=M-\{ o\}$ by $v=ru$.
Fix a center $o=o\0{\HH^k}\in \HH^k\in \cE_k(M)$ and
we get a center $o=o\0{\cE_k(M)}=(o\0{\HH^k},o\0{M})$ of
$\cE_k(M)$.  Then, for $y\in\HH^k-\{ o\}$ we can also write $y=t\,w$, $(w,t)\in \bS^{k-1}\times\R^+$. 
Similarly, using the exponential map we can identify $\cE_k(M)-\{ o\}$
with $\bS^{k+n-1}\times \R^+$, and for $p\in\cE_k(M)-\{ o\}$
we can write $p=s\,x$, $(x,s)\in\bS^{k+n-1}\times\R^+$.\vspace{.1in}

As before denote the metric on $\cE_k(M)$ by $g$ and we can write $g=g_s+ds^2$.
Since $\HH^k$ is convex in $\cE_k(M)$ we can write $\HH^k-\{ o\}=\bS^{k-1}\times \R^+\sbs\bS^{k+n-1}\times \R^+$
and $\bS^{k-1}\sbs \bS^{k+n-1}$.\vspace{.1in}

A point $p\in\cE_k(M)\, -\, o$ has two sets of coordinates: the {\it polar coordinates}
$(x,s)=(x(p),s(p))\in \bS^{k+n-1}\times \R^+$ and the {\it hyperbolic extension coordinates} $(y,v)=(y(p), v(p))\in \HH^k\times M$. Write $M_o=\{o\}\times M$.
Therefore we have the following functions:
$$
\begin{array}{lll}
{\mbox{the distance to {\it o} function:}}  & s:\cE_k(M)\ra [0,\infty), & s(p)=d\0{\cE_k(M)}(p,o)\\\\
{\mbox{the direction of {\it p} function:}}  & x:\cE_k(M)-\{o\}\ra \bS^{n+k-1}, & p=s(p)x(p)\\\\
{\mbox{the distance to {\it $\HH^k$} function:}}  & r:\cE_k(M)\ra [0,\infty), & r(p)=d\0{\cE_k(M)}(p,\HH^k)\\\\
{\mbox{the projection on $\HH^k$ function:}}  & y:\cE_k(M)\ra \HH^k, &\\\\
{\mbox{the projection on $M$ function:}}  & v:\cE_k(M)\ra M, & \\\\
{\mbox{the projection on $\bS^{n-1}$ function:}}  & u:\cE_k(M)-\HH^k\ra \bS^{n-1}, & v(p)=r(p)u(p)\\\\
{\mbox{the length of $y$ function:}}  & t:\cE_k(M)\ra [0,\infty), & t(w)=d_{\HH^k}(y,o)\\\\
{\mbox{the direction of $y$ function:}}  & w:\cE_k(M)-M_o\ra \bS^{k-1}, & y(p)=t(p) w(p)
\end{array}
$$

Note that $r=d_M(v, o)$. Note also that, by 3.1, the functions $w$ and $u$ are constant on geodesics emanating from $o\in\cE_k(M)$, that is
$w(sx)=w(x)$ and $u(sx)=u(x)$.\vspace{.1in}

Let $\p_r$ and $\p_s$ be the gradient vector fields of $r$ and $s$, respectively. Since the $M$-fibers $M_y=\{ y\}\times M$ are convex
the vectors $\p_r$ are the velocity vectors of the speed one geodesics of the form $a\mapsto (y, a\, u)$, $u\in\bS^{n-1}\sbs M$. These geodesics
emanate from (and orthogonally to) $\HH^k\sbs \cE_k(M)$.
Also the vectors  $\p_s$ are the velocity vectors of the speed one geodesics 
emanating from $o\in\cE_k(M)$. For $p\in\cE_k(M)$, denote by $\bigtriangleup =\bigtriangleup (p)$ the right triangle with vertices $o$, $y=y(p)$, $p$
and sides the geodesic segments $[o,p]\sbs\cE_k(M)$, $[o,y]\sbs\HH^k$, $[p,y]\sbs\{ y\}\times M\sbs\cE_k(M)$.
(These geodesic segments are unique and well defined because:\, (1) $\HH^k$ is
convex in $\cE_k(M)$,\, (2) $(y,o)=o_{_{\{ y\}\times M}}$ and $o$ are centers in $\{ y\}\times M$ and $\HH^k\sbs\cE_k(M)$, respectively.)
\vspace{.1in}

\noindent {\bf Lemma 4.1} {\it Let $\eta^+$  (or $\eta$) be a geodesic ray (line) in $M$ through $o$ containing
$v=v(p)$ and $\gamma$ a geodesic line in $\HH^k$ through $o$ containing $y=y(p)$. Then $\bigtriangleup (p)\sbs \gamma\times \eta^+\sbs \gamma\times \eta$.}\vspace{.1in}

\noindent {\bf Proof.} We have that $[o,v]\sbs\eta$ and $[o, y]\in\gamma$. By Lemma 3.1 we have $[o,p]\in \gamma\times \eta^+$.
This proves the lemma.\hfill $\Box$\vspace{.1in}

Let $\alpha:\cE_k(M)-\HH^k\ra \R$ be the angle from 
$\p_s$ to $\p_r$ (in that order), thus  $cos\, \alpha=g(\p_r,\p_s)$, $\alpha\in [0,\pi]$. 
Then $\alpha=\alpha(p)$ is the interior angle, at $p=(y,v)$, of the right triangle $\bigtriangleup =\bigtriangleup (p)$.
We call  $\beta(p)$ the interior angle of this triangle at $o$, that is $\beta(p)=\beta(x)$ is the spherical distance 
between $x\in \bS^{k+n-1}$ and the totally geodesic sub-sphere $\bS^{k-1}$. Alternatively, $\beta$ is the angle between the geodesic segment
$[o,p]\sbs\cE_k(M)$ and the convex submanifold $\HH^k$.
Therefore $\beta$ is constant on geodesics emanating from $o\in\cE_k(M)$, that is
$\beta(sx)=\beta(x)$.\vspace{.1in}

Note that the right geodesic triangle $\bigtriangleup (p)$ has sides of length $r=r(p)$, $t=t(p)$ and $s=s(p)$. By Lemma 4.1 and Remark 3.3
we can consider $\bigtriangleup$ as contained in hyperbolic 2-space.
Hence using hyperbolic trigonometric identities
we can find relations between $r$, $t$, $s$, $\alpha$ and $\beta$. For instance, using the hyperbolic law of cosines we get:\vspace{.1in}

\noindent {\bf (4.2.)}\hspace{1.85in}
$\cosh\, (s)\, =\, \cosh\, (r)\,\, \cosh\, (t)$\vspace{.1in}

\noindent Note that this implies $t\leq s$.
Here is an application of this equation.\vspace{.1in}

\noindent {\bf Proposition 4.3 (Iterated hyperbolic extensions)} {\it We have that $$\cE_l\big( \cE_k(M)   \big)=\cE_{l+k}(M)$$
where we are identifying $\HH^{l+k}$ with $\HH^l\times\HH^k$ with warped metric $(\cosh^2t)\,\sigma\0{\HH^l}+\sigma\0{\HH^k}$.}\vspace{.1in}

\noindent {\bf Remarks.} 

\noindent {\bf 1.} Note that the identification of
$\HH^{l+k}$ with $\HH^l\times\HH^k$ (with warp metric) depends on the order of $l$ and $k$, that is, on the order in which
the hyperbolic extensions are taken.

\noindent {\bf 2.} As before, here the function $t:\HH^k\ra [0,\infty) $ is the distance in $\HH^k$ to the  point $o\in\HH^k$. \vspace{.1in}

\noindent {\bf Proof of Proposition 4.3.} As above let $s:\HH^k\times M\ra [0,\infty)$ be the distance in $\cE_k(M)$ to $o$,
$r(p)=d_M(v(p),o)$, and $t$ as in the statement of the proposition.
Then $\cE_l\big( \cE_k(M)   \big)$ is $\HH^l\times (\HH^k\times M)$ with metric $$(\cosh^2 s)\,\sigma\0{\HH^l}+\big[  (\cosh^2r)\,\sigma\0{\HH^k} +h \big]$$
On the other hand $\cE_{l+k}(M)$ is  $(\HH^l\times \HH^k)\times M$ with metric $$(\cosh^2 r)\,\big[(\cosh^2t)\, \sigma\0{\HH^l}+  \sigma\0{\HH^k}\big] +h $$
\noindent Hence the Proposition follows from identity (4.2) above. This proves the Proposition.\hfill $\Box$\vspace{.1in}

\noindent {\bf Proposition 4.4.} {\it  We have the following identity defined outside $\HH^k\cup \big(\{o \}\times M$\big)}
$$
\big( \sinh^2 (s) \big)\, d\beta\,^2\,\,\, +\,\,\, ds^2\,\,\, =\,\,\, \cosh^2 (r) \, dt^2\,\,\,+\,\,\, dr^2
$$
\noindent {\bf Proof.} First a particular case. Take $M=\R$ and $k=1$, hence $\cE_k(M)=\cE_1(\R)=\HH^2$. In this case the left-hand side of
the identity above is the expression of the metric of $\HH^2$ in polar coordinates $(\beta, s)$, and right hand side of the equation is
the expression of the same metric in the hyperbolic extension coordinates $(r,t)=(v,y)$. (Here $r$ and $t$ are ``signed" distances.)
Hence the equation holds in this particular case.\vspace{.1in}

Now, the general case can be reduced to this particular case using Lemma 3.1 and Remark 1 in 3.3. This proves the proposition.\hfill $\Box$\vspace{.1in}

A direct (and longer) proof of the lemma above can be given using hyperbolic trigonometric identities.
\vspace{.3in}

\noindent {\bf 5. The hyperbolic extension of an $\epsilon$-close to hyperbolic metric.} \\
Let $(M^n,h)$ have center $o$ as above and consider the hyperbolic extension $\cE_k(M)$. As before the metric $\cE_k(h)$ on $\cE_k(M)$ is denoted by $g$.
The ball of radius $a$ on $M$ centered at $o$ will be denoted by $B_a$.
Choose $o\in\HH^k\sbs \cE_k(M)$.
Recall that $o$ is a center of $\cE_k(M)$ (see 3.3), hence we can express $\cE_k(M)-\{ o\}$ as $\bS^{n+k-1}\times\R^+$ with
variable metric $g=g_s+ds^2$.
\vspace{.1in}

Before we prove Theorem A we shall give a particular type of charts on $\HH^{k+1}=\cE_k(\R)$.
For this we consider two ways of describing $\HH^{k+1}$:
using the hyperbolic extension coordinates and
the polar coordinates described in Section 2 (for the case $M=\R$). Hence a point $z\in\HH^{k+1}$ has hyperbolic extension coordinates $(y,r)
\in\HH^k\times\R $ and polar coordinates $(x,t)\in\bS^{k}\times\R^+$.
Using the law of sines and the laws of cosines for right triangles in $\HH^2$ we can find transformation rules between the coordinates $(x,t)$
and the coordinates $(y,r)$. We are only interested in the explicit expression for $r=r(x,t)$. In this case we have\vspace{.1in}

\noindent {\bf (5.1)}\hspace{3.3in}$ r(x,t)\, =\ \sinh^{-1}\bigg(   \sinh\, ( t)\, sin\, \beta (x) \bigg)$\vspace{.1in}

\noindent where $\beta (x)$ is the spherical distance from $x\in \bS^k$ to the equator $\bS^k\cap \HH^k\sbs\HH^{k+1}$.\vspace{.1in}

Fix $\xi\geq 0$. Let $z_0=(x_0,t_0)\in\bS^k\times (2+\xi,\infty)\sbs\cE_k(\R)$ and let $(y_0,r_0)$ be the hyperbolic extension coordinates of $z_0$. The following definition is a particular
version of definition $(\star)$ given at the beginning of the proof
of Theorem 4.2 in \cite{O1}.
We define the chart $\psi=\psi_{z_0}:\T^{k+1}_\xi
\ra\bS^k\times\R^+=\HH^{k+1}-\{ o\}$ by\vspace{.1in}

\noindent {\bf (5.2)}\hspace{1.6in}$ \psi (x,t)\,=\,\bigg( exp_{x_0}\big(e^{\lambda-t_0}\, x\big) \, ,\,  t_0+t \bigg)$\vspace{.1in}

\noindent 
where: (1)  we are identifying the euclidean unit ball $\B^k$ with the unit ball in the tangent space $T_{x_0}\bS^k$, (2)
 $exp_{x_0}:T_{x_0}\bS^k\ra\bS^k$ is the exponential map,
 and (3) $\lambda=min\{0, t_0-ln(k c\0{4})\}$, $c\0{4}=\sqrt{k\,k!\, c^{k}(\bS^k)}$. (The $\lambda$ is a correcting term for $t$ small, see proof
 of 4.2 in \cite{O1}. Here $c(\bS^k)$ is such that $\bS^k$ is 
 $c(\bS^k)$-bounded, see Section 4 of \cite{O1}). Note that in the formula above the output of the map $\psi$ is given in polar coordinates.\vspace{.1in}

\noindent {\bf Lemma 5.3.} {\it The chart $\psi$ is a radially  
$\epsilon$-close to hyperbolic chart, provided}$$C_4\,e^{-t_0}\,\leq\,\epsilon$$
\noindent {\it where $C_4=C_4(k,\xi)$ is as in Remark 4 after Theorem A in the Introduction.}\vspace{.1in}

\noindent {\bf Proof.} 
This lemma was proven in a more general form in the proof of Theorem 4.2 of \cite{O1}: the chart $\psi$ is a special case of the
chart $\phi$ that appears in equation $(\star)$ at the
beginning of the proof of 4.2 \cite{O1}. To see this note: (1) in the proof of
4.2  \cite{O1} it is proven that the chart $\psi$ in $(\star)$  of  \cite{O1} is $\eta$-close to
hyperbolic, provided $\eta\geq C(e^{-t_0}+\epsilon)$, for certain constant $C$ (different from the $C$ in Theorem A),
(2) in our case the map $\varphi$ in $(\star)$ is the exponential map
$exp_{x_0}$, hence
$A$ in $(\star)$ is the identity (the derivative of the exponential is the identity),
(3) we can take $\epsilon=0$ in $(\star)$
because the family of metrics $\{\sigma\0{\bS^k}\}$ is constant
hence ``zero"-slow (i.e $\epsilon$-slow for every $\epsilon>0$, see Section 4 of \cite{O1}), (4) the chart in $(\star)$ works to prove
Theorem 4.2 \cite{O1} (where warping function $e^t$ is used), and the same chart works to prove Corollary 4.3 \cite{O1}
(where warping function $\sinh t$ is used), but the constant
$C$ changes to a new constant $C_1=C_1(c,k,\xi)$ (see item 3 in Remarks
1.3 \cite{O1}, Remark 3 after the statement of Theorem 4.2 \cite{O1} and
Remark 2 after the statement Corollary 4.3 \cite{O1}).  In our case
we can take $c=c\0{\bS^k}$ such that $\bS^k$ is $c\0{\bS^k}$-bounded, therefore our constant becomes $C_4(k,\xi)=C_1(c\0{\bS^k},
k,\xi)$ as in Remark 4 after Theorem A in the Introduction.
This proves the lemma. \hfill $\Box$
\vspace{.1in}

Denote the hyperbolic extension coordinates of \, $\psi$ by \,$y=y_{z_0}:\T_\xi\ra \HH^k$ and \,$\br=\br_{z_0}:\T_\xi\ra \R$. That is
$$\psi (x,t)\,\,=\,\,\big( y(x,t)\, ,\, \br(x,t)\big)\, \in\, \HH^k\times \R\, =\, \cE_k(\R)
$$
\noindent Using equation (5.1) we can write\vspace{.1in}

\noindent {\bf (5.4)}\hspace{3.1in} $\br(x,t)\, =\ \sinh^{-1}\bigg(   \sinh\, ( t_0+t)\, sin\, \beta (x') \bigg)$\vspace{.1in}

\noindent where $x'=exp_{x_0}\big(e^{\lambda-t_0}\, x\big) $. 
Recall that $(y_0,r\0{0})$ are the hyperbolic extension coordinates of $z_0$.\vspace{.1in}

\noindent {\bf Lemma 5.5.} {\it We have that}
$$\big|\, \br (x,t)\,\, -\,\, \big(  t\, +\, r\0{0} \big)\,\big|_{C^2}\,\, \leq \,\, Le^{-r\0{0}}$$
\noindent {\it where $L$ is a constant depending on $k$ and $\xi$.}\vspace{.1in}

The proof of Lemma 5.5 is given in Appendix A.
An explicit formula for $L$ is given at the end of Appendix
A.\vspace{.1in}

The next result is the reason why we introduced the variable $\xi$ in 
the definition of the models $\T_\xi$: the new excess $\xi'$  
is less than the old excess $\xi$.
\vspace{.1in}

\noindent {\bf Corollary 5.6.} {\it We have that}
$$ \psi (\T_{\xi'})\,\sbs\, \HH^k\times \big[r\0{0}-(1+\xi)\, ,\, r\0{0}+(1+\xi)\big ] 
$$\noindent {\it provided $0<\xi' <\xi-
Le^{-r\0{0}}$, where $L$ is as in 5.5}.\vspace{.1in}

\noindent {\bf Proof.} Write $\kappa =
Le^{-r\0{0}}$. By Lemma 5.5 we have $(t+r\0{0})-\kappa\leq \br(x,t)\leq (t+r\0{0})+\kappa$.
Hence for $t\in (-1-\xi', 1+\xi')$ we get $r\0{0}-(1+\xi'+\kappa)\leq \br(x,t)\leq r\0{0}+(1+\xi'+\kappa)$. This together with the choice  $\xi'+\kappa\leq \xi$ implies $r\0{0}-(1+\xi)\leq \br(x,t)\leq r\0{0}+(1+\xi)$. This proves the corollary.\hfill $\Box$\vspace{.1in}

\noindent {\bf Proof of Theorem A.} First some notation. 
Recall that we are denoting the metric on $M-\{ o\}$ by $h=h_r+dr^2$ and the one on $\cE_k(M)$ by $g$.
Also recall $I_\xi =[-1-\xi,1+\xi]$. For $u\in \bS^{n-1}$ we denote by
$Ru$ the complete geodesic line in $M$ passing through $o$ with direction $u$, i.e. $Ru=exp_o (\R u)=\{ p\in M$ such that $p=ru$, i.e. 
$p$ has polar coordinates $(u,r),\, r\in\R \}$. Also write $R^+u=\exp_o(\R^+u)$. Then $Ru=R(-u)$ but $R^+u\cap R^+(-u)=\emptyset$.
Hence we get  $M-\{ o\}=\coprod _{u\in\bS^{n-1}} R^+u$. Therefore
\begin{equation*} 
\cE_k(M)-\HH^k\,\,\, =\,\,\, \HH^k\, \times \, \big(M-\{ o\}\big)\, \,\,=\, \,\,\coprod_{u\in\bS^{n-1}}\bigg( \HH^k\times R^+u\bigg)
\tag{1}\end{equation*}
\noindent Specifically, if $w=(y,p)\in \HH^k\times M=\cE_k(M)$, $p\neq o$, then $w\in \HH^k\times R^+u$, provided $p$ has polar coordinates $(u,r)$, for some $r>0$.\vspace{.1in}

Since we are taking $u$ with length one,
we have an obvious identification of $Ru$ with $\R$, given by $r\mapsto ru$ (this identification does depend on the ``sign" of $u$).
This identification gives a canonical (metric) identification of $\cE_k(Ru)=\HH^k\times Ru$ (with metric $g|_{\cE_k(Ru)}$) 
\, with\, $\cE_k(\R)=\HH^k\times\R=\HH^{k+1}$ (with the canonical warped product metric).\vspace{.1in}

Write $\cE_k(R^+u)=\HH^k\times R^+u\sbs\cE_k(Ru)$ and
note that we can canonically identify $\cE_k( R^+u)$\, with half hyperbolic $(k+1)$-space \, 
 $\HH^{k+1}_+=\HH^k\times \R^+\sbs \HH^k\times \R=\HH^{k+1}$.
\vspace{.1in}

For $r>0$ and $y\in \HH^k$ denote by $\bS_{r,y}$ the set $\{ (y,ru)\in \HH^k\times M,\,\,\, u\in\bS^{n-1}\,  \}$. Then $\bS_{r,y}$ is the geodesic
sphere of radius $r$ of the convex submanifold $\{ y\}\times M\sbs\cE_k(M)$. Note that every $\cE_k(R^+u)$ intersects every $\bS_{r,y}$ orthogonally in the single point $(y, ru)$.\vspace{.1in}

Let $w_0\in \HH^k\times S \sbs\cE_k(M)$. Write $w_0=(y_0,p_0)\in \HH^k\times M$ and let $(u_0,r\0{0})$ be the polar coordinates of $p_0\in S\sbs M$.
Also let $t_0$ be the distance in $\cE_k(M)$ from $o$ to $w_0$.
Since $p_0\in S$ and $S$ is radially $\epsilon$-close to hyperbolic, there is an radially $\epsilon$-close to hyperbolic chart  $\phi : \T^{n}_\xi\ra M$ with center $p_0$. From the definition of a radially $\epsilon$-close to hyperbolic chart (see Section 3) we have that 
for $(x,r)\in \T^n_\xi=\B^{n-1}\times I_\xi$ we can write 
\begin{equation*}
\phi\,(x,r)\,\,=\,\,\big(u(x), r+r\0{0}\big)=(r+r\0{0})\,u(x)
\tag{2}
\end{equation*}
Also write $z_0=( y_0, r\0{0})\in\cE_k(\R^+u)\sbs\cE_k(\R u)
=\HH^{k+1}$ and let $\psi=\psi_{z_0}:\B^k\times\R^+\ra\HH^{k+1}$, be
as defined in (5.2).
And, as before, write $\psi=\psi_{z_0}=(y_{z_0}, \br_{z_0})=(y,\br)$,
where $(y,\br)$ are the hyperbolic extension coordinates of $\psi$. Note that we are taking the domain of $\psi $ as $\B^k\times\R^+$ and not just $\T_\xi^{k+1}=
\B^k\times I_\xi$ (see Remark 2 at the end of Section 3).
We now define a chart $\bphi:\B^k\times\B^{n-1}\times \R^+\ra \cE_k(M)$ by
\begin{equation*}
\bphi \,(x_1, x_2, t)\,\,=\,\,\bigg(y(x_1,t)\, ,\,  \br(x_1,t)\,u(x_2)\bigg)\, \in \HH^k\times M\,=\,\cE_k(M)
\tag{3}
\end{equation*}
\noindent Note that, by Corollary 5.6 we have that 
\begin{equation*}
\bphi\,\big(    \B^k\times\B^{n-1}\times I_{\xi'}     \big)\, \sbs\, \HH^k\times \phi (\T_\xi)
\tag{4}
\end{equation*}
\noindent provided $\xi'<\xi-Le^{r\0{0}}$. By the definition of $\bphi$ (see equation (3)) and (2) we have
\begin{equation*}
\bphi \bigg(\{ x_1\}\times\B^{n-1}\times \{ t\}\bigg)\,\,=\,\, \phi \bigg( \B^{n-1}\times \{\br(x_1,t)-r\0{0} \}  \bigg)\,\,\sbs\,\,
\bS_{\br(x_1,t), y(x_1,t)}
\tag{5}
\end{equation*}
\noindent and 
\begin{equation*}
\bphi \bigg(\B^{k}\times \{ x_2\}\times I_{\xi'}\bigg)\,\,\sbs\,\,\cE_k\big(    R^+u(x_2) \big)
\tag{6}
\end{equation*}\vspace{.1in}

\noindent Moreover, using (6), the canonical metric identification 
between  $\cE_k(Ru)$ and $\cE_k(\R)=\HH^{k+1}$,
and the obvious identification between $\B^{k}\times \{ x_2\}\times I_{\xi'}$ and $\T_{\xi'}= \B^{k}\times I_{\xi'}$,
we can say that the chart $\bphi$ satisfies
\begin{equation*}
\bphi \big| _{\B^{k}\times \{ x_2\}\times I_{\xi'}}\,\,=\,\, \psi 
\tag{7}
\end{equation*}
\noindent Also, from (5), (3), (2), and using the obvious identifications of $\{ x_1\}\times \B^{n-1}\times\{ t\}$\, with $\B^{n-1}\times\{ \br (x_1,t)-r\0{0}\}$
\, and $\{ y(x_1,t)\}\times M$\, with  $M$, we can write
\begin{equation*}
\bphi \big| _{\{ x_1\}\times \B^{n-1}\times\{ t\}}\,\,=\,\, \phi\big|_{\B^{n-1}\times  \{  \br(x_1,t)-  r\0{0}  \} }
\tag{8}
\end{equation*}
Since, as mentioned above, every $\cE_k(R^+u)$ intersects 
every $\bS_{r,y}$\, $g$-orthogonally in a single point, we have that the $\B^{n-1}$-fibers \,\, $\{ x_1\}\times\B^{n-1}\times \{ t\}$\,,  and the \,$(\B^k\times I)$-fibers
\,\,$\B^{k}\times \{ x_2\} \times I_{\xi'}$\,\, are $\bphi^*(g)$-orthogonal. Also, by (7),  $\bphi^*(g)$ restricted to a $(\B^k\times I)$-fiber is canonically
hyperbolic, hence, by Lemma 5.3 and the fact that $r\0{0}\leq t_0$, we have\vspace{.1in}

\noindent {\bf (5.7)} \,\,\,{\it
The map $\bphi$, restricted to a $\B^k\times I$-fiber, 
 is a $\delta$-close to hyperbolic chart, provided $C_4\,e^{-r\0{0}}\leq \delta$}\vspace{.1in}

Therefore $\bphi$ has the form $\bphi^*g=f_1+f_2+dt^2$
where $f_1$ is the restriction of $\bphi^*g$ to the $\B^k$-fibers and $f_2$ is the restriction of $\bphi^*g$ to the $\B^{n-1}$-fibers.
Also, $f_1+dt^2$ is the restriction of $\bphi^*g$ to the $(\B^k\times I)$-fibers, and $f_1+dt^2$ is  a hyperbolic metric.
Furthermore, again by (7), we have that  $f_1+dt^2$ (hence also $f_1$) 
is independent of the variable $x_2$. Now we need an estimate
for $f_2$.\vspace{.1in}

\noindent {\bf Claim 5.8.} {\it We have that}
$$\big|\, f_2(x_1,x_2,t)\,\, -\,\,e^t\sigma\0{\R^{n-1}} \big|_{C^2}\,\, \leq\,\, 2(2+3\xi+\xi^2)e^{^{1+\xi}}L\,\Big( \epsilon\,+\, e^{-r\0{0}}\Big)$$\vspace{.1in}

\noindent {\bf Proof of Claim 5.8.} Let $a_{ij}$ be the entries of the matrix $f_2$. We have to prove that
$\big|\, a_{ij}(x_1,x_2,t)\,\, -\,\,e^t\delta_{ij} \big|_{C^2}\,\, <\,\, 
2(2+3\xi+\xi^2)e^{^{1+\xi}}L\,\big( \epsilon\,+\, e^{-r\0{0}}\big)$.
Let $b_{ij}(x,r)$ be the entries of the matrix $\phi^*h_r$. Since, by hypothesis, $\phi$ is an $\epsilon$-close to hyperbolic
chart, we
have that 
\begin{equation*}
\big|\, b_{ij}(x\,,\,r)\,\, -\,\,e^r\delta_{ij} \big|_{C^2}\,\, <\,\, \epsilon 
\tag{9}
\end{equation*}
\noindent On the other hand, equation (8) implies:
\begin{equation*}
a_{ij}(x_1\,,\,x_2\,,\,t)\,=\, b_{ij}\big(x_2\, ,\,\br(x_1\, ,\,t)\,-\, r\0{0}\,\big)
\tag{10}
\end{equation*}
The proof of the claim is obtained by calculating the derivatives of $a_{ij}(x_1,x_2,t)\,\, -\,\,e^t\delta_{ij}$ 
up to order 2 and finding estimates of these derivatives using (9), (10) and Lemma 5.5.
This is done in Appendix B. (The idea here is that, by (10), $a_{ij}(x_1,x_2,t)\,\, -\,\,e^t\delta_{ij}$ 
is equal to $b_{ij}(x_2,\br (x_1,t)-r\0{0})\,\, -\,\,e^t\delta_{ij}$, which, by Lemma 5.5 is $C^2$-close to
$b_{ij}(x_2,t)\,\, -\,\,e^t\delta_{ij}$ which, by (9), is small.)\hfill $\Box$\vspace{.1in}

We now complete the proof of Theorem A. 
Choose $\delta=C_4e^{-r\0{0}}$ in 5.7.
Recall that  $f_1+dt^2$ is the restriction of $\bphi^*g$ to the $(\B^k\times I)$-fibers. This together with
5.7 and 5.8 imply the following (all norms are $C^2$)
$$\begin{array}{lll}\big|\,\bphi^*g-\big(e^t\sigma\0{\R^{n+k-1}}+dt^2\big)\,\big|&=&\big|\,f_1+f_2+dt^2-\big(e^t\sigma\0{\R^{k}}+
e^t\sigma\0{\R^{n-1}}+dt^2\big)\,\big|\\\\
&\leq&\big|\,f_1+dt^2-\big(e^t\sigma\0{\R^{k}}+dt^2\big)\,\big|
\,+\,\big|\,f_2-e^t\sigma\0{\R^{n-1}}\,\big|\\\\
&\leq&\, C_4\,e^{-r\0{0}}\,\,\,+\,\,\,2(2+3\xi+\xi^2)e^{^{1+\xi}}L\,\Big( \epsilon\,+\, e^{-r\0{0}}\Big)\\\\
&\leq&\big[2(2+3\xi+\xi^2)e^{^{1+\xi}}L\,+\, C_4\big]\,\,\Big( \epsilon\,+\, e^{-r\0{0}}\Big)
\end{array}$$

\noindent Finally recall that we had chosen $w_0\in \cE_k(S) \sbs\cE_k(M)$, with $w_0=(y_0,p_0)\in \HH^k\times M$ and $(u_0,r\0{0})$ the polar coordinates of $p_0\in S\sbs M$. If $w_0\notin \cE_k(B_a)\sbs\cE_k(M)$, then $r\0{0}\geq a$, thus $e^{-r\0{0}}\leq e^{-a}$.
This completes the proof of Theorem A.\hfill $\Box$\vspace{.1in}

\noindent {\bf Proposition 5.9.}  {\it  Let $M^n$ have center $o$. Assume $M$ is $(B_a,\epsilon)$-close to hyperbolic, with charts of excess $\xi>0$. Then $\cE_k(M)$ is $(B_a,\eta)$-close to hyperbolic, with charts of excess $\xi'$,
provided

$$
C_1'\,e^{1+\xi}\, e^{-a}\,\,+\,\,C\,\, \epsilon\leq\eta
$$
\noindent where $0<\xi'<\xi-Le^{1+\xi}e^{-a}$. Here $C_1'=C_1'(n+k,\xi)$, and $C$, $L$ are as in Theorem A.}\vspace{.1in}

\noindent {\bf Remark.}
 Here $C_1'=C_1'(n+k, \xi)$ is as in 4.13 of \cite{O1}.
The space $\cE_k(M)$ in the proposition above is radially $\eta$-close to hyperbolic
with respect to any center $o\in\HH^k\sbs\cE_k(M)$.
\vspace{.1in}

\noindent {\bf Proof.} Denote the center of $\cE_k(M)$ by $o=(o\0{\HH^k},o\0{M})$. Recall that $\cE_k(\qH^n)=\qH^{n+k}$.
Since $M$ is (radially) hyperbolic on $B_a$
we have
\begin{equation*}
{\mbox{\it  the space $\cE_k(M)$ is (radially) hyperbolic on 
$\HH^k\times\B_{a}(M)=\cE_k(B_a)$}}
\tag{a}
\end{equation*}
\noindent Let $p=(y,v)\in\cE_k(M)$. We use the functions (coordinates) in Section 2. 
In particular $s=d\0{\cE_k(M)}(o,p)$, 
$r=d\0{\cE_k(M)}(p,\HH^k)=d\0{M}(o\0{M},v)$.  
Recall that $s\geq r$ (see 4.2). 
We will use Corollary 4.13 of \cite{O1}: \vspace{.1in}

\noindent {\it Corollary 4.13 of \cite{O1}. There is
 $C_1'=C_1'(n+k,\xi)$ such that
hyperbolic $(n+k)$-space $\qH^{n+k}$ is
radially $(C_1'e^{-b})$-close to hyperbolic outside $B_b$.}
\vspace{.1in}

\noindent We have two cases.\vspace{.1in}

\noindent {\bf First case.} $r\geq a-1-\xi$.

\noindent It can be checked that $C_1'(n+k,\xi)>C(k,\xi)$,
hence the hypothesis
$C_1'e^{1+\xi}\,\, e^{-a}\,\,+\,\,C\,\, \epsilon\leq\eta$ implies
$C\,\big( e^{-(a-1-\xi)}\,\,+\,\, \epsilon\big)\leq\eta$. This together with
Theorem A (replace $a$ by $a-1-\xi$) imply that $\cE_k(M)$ is radially $\eta$-close to hyperbolic at $p$
(i.e there is a radially $\eta$-close to hyperbolic chart centered at $p$, with excess $\xi'<\xi-Le^{1+\xi}e^{-a}$).
\vspace{.1in}

\noindent {\bf Second case.} $r < a-1-\xi$, $s\geq a-1-\xi$.\\
The hypothesis
$C_1'\,e^{1+\xi}\, e^{-a}\,\,+\,\,C_3\,\, \epsilon\leq\eta$ implies
$C_1'\,e^{1+\xi} e^{-a}\leq\eta$.  Since $s\geq a-1-\xi$ we can  apply
Corollary 4.13 of \cite{O1} (stated above; take $b=a-1-\xi$) to obtain that $\qH^{n+k}$ is radially $\eta$-close to hyperbolic outside the ball $\B_{a-1-\xi}(\qH^{n+k})$ of radius $a-1-\xi$
on $\qH^{n+k}=\cE_k(\qH^n)$. By (a) we can identify
$\cE_k(B_a)\sbs\cE_k(M)$ with the corresponding subset
of $\qH^{n+k}=\cE_k(\qH^n)$. Since $p\in\cE_k(B_{a-1-\xi})\sbs\cE_k(M)$, we can also consider $p\in\cE_k(B_{a-1-\xi})\sbs\cE_k(\qH^n)$. Therefore there is a radially $\eta$-close to hyperbolic chart $\phi$ centered at $p$,
with image (a priori) contained in $(n+k)$-hyperbolic space.
Note that, since $\phi$ is centered at $p$ and we are assuming $r< a-1-\xi$, by Corollary 5.6 we get that
$\phi(\T_{\xi'})\sbs \{ r< a\}$.
Hence, again by (a), the chart $\phi$ is also a chart for $\cE_k(M)$.
This proves the proposition.\hfill $\Box$\vspace{.1in}

We now prove Theorem B given in the Introduction.\vspace{.1in}

\noindent {\bf Proof of Theorem B.} 
We have $a\geq R=ln\big(\frac{1}{\epsilon}\big)+ln(L)+1+\xi$. 
Therefore we get:
(1) $e^{-a}\leq \epsilon$, and (2) $e^{a/2}\leq Le^{1+\xi}$.
Note that (2) implies $\xi-e^{-a/2}\leq \xi-Le^{1+\xi}e^{-a}$.
We can now apply 5.7. Note that we are writing
$C_2=C_2(n,k,\xi)=C_1'e^{1+\xi}+C$.
This proves the Theorem.\hfill $\Box$
\vspace{.2in}

\noindent {\bf \large  Appendix A. Proof of Lemma 5.5.}\\
Recall that we are considering $\HH^{k+1}$ with two sets of coordinates: the polar coordinates
$(x,t)$ and the hyperbolic extension coordinates $(y,r)$. Recall that $x,t,y,r$ are functions defined on $\HH^{k+1}$, specifically:
$x:\HH^{k+1}-\{ o\}\ra \bS^k$, $y:\HH^{k+1}\ra\HH^k$,   $r:\HH^{k+1}\ra\R$,    $t:\HH^{k+1}\ra\R$.
Let $\p_r$ be the gradient vector field of $r$. Then the vectors $\p_r$ are the velocity vectors of the speed one geodesics 
emanating orthogonally from
$\HH^k\sbs \HH^{k+1}$. Also let $\p_t$ be the gradient vector field of $t$ and let $\alpha:\HH^{k+1}-\{ o\}\ra \R$ be the angle between 
$\p_t$ and $\p_r$. Then $\alpha(z)$ is the interior angle, at $z=(y,r)$, of the right triangle with vertices $o$, $y$, $z$. We call 
$\beta(z)$ the interior angle of this triangle at $o$, that is $\beta(z)=\beta(x)$ is the  (signed) spherical distance 
between $x\in \bS^k$ and the equator $\bS^{k-1}\sbs\bS^k$, where $(x,t)$ are the polar coordinates of $z$.
Note that the triangle mentioned above has sides of length $r=r(z)$, $t=t(z)$ and $a=a(y)$, where we are denoting by $a$ the distance
function in $\HH^k$ to $o$.
Using the hyperbolic law of cosines we get:
\begin{equation*}
sin\,\alpha\, =\, \frac{cos\, \beta}{\cosh\, r}
\tag{A1}
\end{equation*} 
\noindent Therefore 
\begin{equation*}
\begin{array}{ccccc}
|sin\,\alpha |\leq \frac{1}{\cosh\, r} &\,\,\,\,\,\, &&{\mbox{and}}&\,\,\,\,\,\,\,\,\,\,\,\,\,\,\,\,
|cos\, \alpha |\geq \frac{|\sinh\, r|}{\cosh\, r}\,\, =\,\,| tanh\, r|
\end{array}
\tag{A2}
\end{equation*} 
\noindent Note that the map $sin\, \beta$ is just the height function, i.e. $sin\,\beta (x)$ is the (signed) euclidean distance from 
$x\in\bS^k$ to $\HH^k$, which is  the last coordinate $x_{k+1}$ of $x=(x_1,...,x_{k+1})$. Therefore the term $sin\, \beta(x')$ that appears in
the definition of $\br$ (see 5.4) is the composition 
$$\B^k\stackrel{ {\tiny e^{\lambda-t_0}}}{\longrightarrow}\,\,\B^k\,\,
\stackrel{exp}{\longrightarrow}\,\,\bS^k\,\,\stackrel{proj}{\longrightarrow}\,\,\R$$
\noindent where the first arrow is multiplication by the constant $e^{\lambda-t_0}$, $exp=exp_{x_0}$, and the last arrow is the projection take-the-last-coordinate map.\vspace{.1in}

Write \,$\p_t=\pt$\,
and \,$\p_i=\frac{\p}{\p x_i}$, $i=1,...,k$. 
Since $exp_{x_0}$ and $proj$ are smooth and the sphere is compact there is a constant $c$ (independent to $x_{0}$)
with $|\, proj\,\circ exp\, |_{C^2}\, \leq c$.\vspace{.1in}

\noindent {\bf Remarks.}

\noindent {\bf 1.} The map $\beta$ is continuous but not smooth at the north pole. On the other hand  $sin\, \beta$ is smooth.

\noindent {\bf 2.} In what follows we will use the fact that we can take $c=3$. 
Moreover we can take $|\nabla(\, proj\,\circ exp\,) |\, \leq 3$
(here $|.|$ is Euclidean length).
A straightforward calculation 
(not given here) can show this.
\vspace{.1in}

\noindent Write $\Lambda=k^{3/2}\sqrt{k!}\,c\0{\bS^k}^{k/2}$, thus $e^\lambda\leq \Lambda$ (actually $e^\lambda=\Lambda$  for
$t_0$ small and $e^\lambda= 1$ otherwise). We have then

\begin{equation*}
\begin{array}{ccccc}
\big|\,\, sin\, \big( \beta \, ( x')\big)\,\,\big|_{C^1}\,\leq\,3\Lambda\,
e^{- t_0}&&\,\,\,\,{\mbox{and}}\,\,\,&&
\big|\,\, sin\, \big( \beta \, ( x')\big)\,\,\big|_{C^2}\,\,\leq 3\Lambda^2\,
e^{-2 t_0}
\end{array}
\tag{A3}
\end{equation*} \vspace{.0in}

\noindent Similarly we have $\big|\nabla sin\, \big( \beta \, ( x')\big)\big|\leq3\Lambda\,
e^{- t_0}$.  Write $\bt=t_0+t$ and note that $\bt>0$ (recall we are assuming $t_0>1+\xi$).
Differentiating equation (5.4) we get\vspace{.1in}

{\small  {\begin{center} $ \bigg| \p_t \br(x,t)\bigg|\, =\, \bigg| \frac{\cosh (\bt)\, sin\big( \beta(x')  \big)}{\cosh\,\br}\bigg|\,\,
=\, \bigg| \frac{\cosh (\br)\, \cosh (a)\, sin\big( \beta(x')  \big)}{\cosh\,\br}\bigg|\,\,
=\,\, cos\,\alpha \,\,\geq\,|\, tanh\, \br \,| $\end{center}}}

\noindent where the second equality  is obtained from the first hyperbolic law of cosines and the last from the second hyperbolic law of
cosines, and the last inequality comes from  (A2). Note also that we get $|\p_t \br|\leq 1$.
Similarly, using further differentiation, (A2), the two laws of cosines, the law of sines, and a bit of work show
{\small {\begin{center}   $ \bigg| \p_t^2 \br(x,t)\bigg|\, =
\, \bigg| \big( tanh\, \br\, \big)\, \big( sin^2\, \alpha\,\big)\bigg| \,\,\leq\,\, \frac{1}{\cosh^2\, \br}  $\end{center}}}
\noindent Also, using (A3) we get
{\small {\begin{center}   $\bigg| \p_i \br(x,t)\bigg|\, =\, \bigg| \frac{\sinh (\bt)\,\p_i sin\big( \beta(x')  \big)}{\cosh\,\br}\bigg|\,\,
\leq \, 3\Lambda \bigg(\frac{\sinh\,(t+t_0)}{e^{t_0}}\bigg)\,\,\frac{1}{\cosh\, \br}\,\,
\leq \,\,3\Lambda (\frac{e^{1+\xi}}{2})\frac{1}{\cosh\, \br}\,\,  $\end{center}}}
\noindent (recall $t\in (-(1+\xi),1+\xi)$). A similar argument
using $\big|\nabla sin\, \big( \beta \, ( x')\big)\big|\leq3\Lambda\,
e^{- t_0}$ shows
{\small {\begin{center}   $\big| \nabla \br(x,t)\big|\,\leq \,\,3\Lambda (\frac{e^{1+\xi}}{2})\frac{1}{\cosh\, \br}\,\,  $\end{center}}}
\noindent Differentiating again,
using the two laws of cosines, and (A3) we obtain
{\small {\begin{center}   $\bigg|\p_t\p_i \br(x,t)\bigg|\,
\leq \,3\Lambda \bigg[\bigg(\frac{\cosh\,(t+t_0)}{e^{t_0}}\bigg)
\,\,+\,\,\bigg(\frac{\sinh\,(t+t_0)}{e^{t_0}}\bigg)\,\bigg]\,\,\frac{1}{\cosh\, \br}\,\,
\leq \,\,\frac{6\Lambda e^{^{1+\xi}}}{\cosh\, \br}\,\,  $\end{center}}}
\noindent provided $t_0\geq 1+\xi$. Finally, differentiating and using
(A3) we get
{\small {\begin{center}   $ \bigg| \p_{ij} \br(x,t)\bigg|\,
\leq \, \Lambda^2\bigg(\frac{\sinh\,(t+t_0)}{e^{t_0}}\bigg)\bigg[ \frac{3}{e^{t_0}}\,\frac{1}{\cosh\, \br}\,\,+\,\,
\frac{\sinh\,(t+t_0)}{e^{t_0}}\frac{9}{\cosh^2\,\br }\,\bigg]\,\,
\leq \,8\Lambda^2\frac{e^{^{2(1+\xi)}}}{\cosh\, \br}$\end{center}}}
\noindent  Note that all five terms on the right
of the last five equations are less than $\frac{8\Lambda^2\,e^{^{2(1+\xi)}}}{\cosh\,\br}\leq16\, \Lambda^2 e^{2(1+\xi)}e^{-\br}$.\vspace{.1in}

Now, write $F(x,t)=\br (x,t)-(t+r\0{0})$. 
Since all but one of derivatives
of order 1 and 2 of $F$ concide with the ones of $\br$ we get that all such
derivatives are less than $16\, \Lambda^2 e^{2(1+\xi)}e^{-\br}$. The remaining derivative is
$\p_t F$. We have
$$|\p_t F|\leq |1-tanh\,\br|=\frac{e^{-\br}}{\cosh\,\br}\leq \frac{1}{\cosh\,\br}\leq16\, \Lambda^2 e^{2(1+\xi)}e^{-\br}$$
It remains to estimate $|F|_{C^0}$.
For $x\in\B^k$, since $F(0,0)=0$, we have $F(x,0)=\int_0^1 x\,.\, \p _x F(tx,0)\, dt$. But $\p_xF=\langle x,\nabla\br\rangle$, hence
$|F(x,0)|\leq \frac{3\Lambda e^{1+\xi}}{2\, \cosh\, (\br)}$.
Hence $$\Big|\,F(x,t)\,\Big|\,=\,\Big|\,F(x,0)+\int_0^1 \p _t F(x,t)\, dt\,\Big|\,\leq\, \frac{3\Lambda e^{^{1+\xi}}}{2\,\cosh\, (\br)}+\frac{1}{\cosh\, (\br)}\,<\, 4\Lambda
e^{^{1+\xi}}e^{-\br}$$
Therefore
\begin{equation*}\big|\br(x,t)-(t+r\0{0})\big|_{C^0}\,\leq\,4\Lambda
e^{^{1+\xi} }e^{-\br}\tag{A4}\end{equation*} \vspace{.1in}

\noindent Hence

\begin{equation*}\big|\br(x,t)-(t+r\0{0})\big|_{C^2}\leq 16\, \Lambda^2 e^{2(1+\xi)}e^{-\br}\tag{A5}\end{equation*} \vspace{.1in}

\noindent To finish the proof we need compare $\br$ with 
$r\0{0}$. For this note that, since $r\0{0}=\br(0,0)$, we have
{\small $$\begin{array}{lll}
\Big| \, \br(x,t)\,-r\0{0}  \, \Big|&\leq&\Big| \, \br(x,t)\,-\,\br(x,0)  \, \Big|
\,+\,\Big| \, \br(x,0)\,-\,\br(0,0)  \, \Big|\\\\
&\leq&\int^t_0\big| \p\0{t}\br(x,t) \big|\, dt\,\,+\,\,
\int^1_0|x|\,\big| \nabla\br(tx,0)\big|\,dt\\\\
&<& (1+\xi)\,+\,\frac{3\Lambda \,e^{^{1+\xi}}}{2\cosh\,\br}\,\leq
(1+\xi)+3\Lambda e^{1+\xi}
\end{array}$$}

\noindent That is

\begin{equation*}
\Big| \, \br(x,t)\,-r\0{0}  \, \Big|_{C^0}\leq
(1+\xi)+3\Lambda e^{1+\xi}\tag{A6}\end{equation*}

\noindent Hence $$\br\,>\,
 r\0{0}\,-\, (1+\xi+3\Lambda e^{^{1+\xi}})$$

\noindent This together with (A5) imply
$$
\big|\br(x,t)-(t+r\0{0})\big|_{C^2}\leq 16\, \Lambda^2 e^{1+\xi}
e^{^{(1+\xi)+3\Lambda e^{1+\xi}}}e^{-r\0{0}}\,=\, Le^{-r\0{0}}
$$

\noindent with $L=L(k,\xi)=16\, \Lambda^2 e^{1+\xi}
e^{^{(1+\xi)+3\Lambda e^{1+\xi}}}$, where $\Lambda^2=k^3k!c^k\0{\bS^k}$.
This completes the proof of the lemma. \hfill $\Box$


\vspace{.2in}

\noindent {\bf \large  Appendix B. Calculations for the proof of Claim 5.8.}\\
We will use the following abbreviations for the partial derivatives: $\p_t=\pt$,
 $\p_i=\frac{\p}{\p u_i}$, $\barp_i=\frac{\p}{\p v_i}$,
where $x_1=(u_1,...,u_k)$ and $x_2=(v_1,...,v_{n-1})$.\vspace{.1in}

Write $\kappa=\kappa(r\0{0})=
Le^{-r\0{0}}$. Note that we are assuming $\kappa\leq \xi$
(see Remark 3 after the statement of Theorem A in the Introduction).
Also write
\begin{equation*}  
\zeta(x_1,t)=\br(x_1,t)-r\0{0}\tag{B1}
\end{equation*}
\begin{equation*}
c=c_{ij}(x_1,x_2,t)= a_{ij}(x_1, x_2,t)-e^t\delta_{ij}=b_{ij}(x_2,\zeta)-e^t\delta_{ij}\tag{B2}
\end{equation*}
\noindent where the last equality follows from (10) in the proof
of Claim 5.8 in Section 4. Also write
\begin{equation*}
d=d(x_2,r)\, =\, b_{ij}(x_2,r)\,-\, e^r\delta_{ij}
\tag{B3}
\end{equation*}
We have to prove that $|c|_{C^2}<2(2+3\xi+\xi^2)e^{^{1+\xi}}L(
\epsilon+e^{-r\0{0}})$.
It follows from (B2), (B3) that
\begin{equation*}
c=d(x_2,\zeta)+e^\zeta\delta_{ij}-e^t\delta_{ij}
\tag{B4}
\end{equation*}

\noindent From Lemma 5.5 and  (9) in the proof of 5.8 in Section 4, we have 
\begin{equation*}
\begin{array}{ccccc}
|\,\zeta\,\, -\,\, t\,|_{C^2}\,\leq \,\kappa &&{\mbox{and}}\,\,\,&& |d|_{C^2}\, <\, \epsilon\,\,\,\,
\end{array}
\tag{B5}
\end{equation*} 
\noindent  It follows from (B1),  Corollary 5.6, and 
$t\in I_{\xi'}$ that
\begin{equation*}
|\zeta|_{C^0}\leq1+\xi
\tag{B6}
\end{equation*}
\noindent Note that (B5) also implies
\begin{equation*}
\begin{array}{cccccccc}
\big| \p_i\zeta \big|_{C^0}\leq\kappa\,\,\,\,\,\,&&
\big| \p_t\zeta \big|_{C^0}\leq 1+\kappa\,\,\,\,\,&&
\big| \p^2_t\zeta \big|_{C^0}\leq\kappa\,\,\,\,\,&&
\barp_j\zeta=0
\end{array}
\tag{B7}
\end{equation*}
\noindent From (B3) and (B5)
we get that for $r\in I_\xi$ we have
\begin{equation*}
\big| \p_r b_{ij}   \big|_{C^1}\,\leq\, \epsilon\,+\, \big| e^r  \big|_{C^2 }\, <\, \epsilon\,  +\, e^{^{1+\xi}}
\tag{B8}
\end{equation*} 
\vspace{.1in}

\noindent {\bf The $C^0$-norm of $c$.}
Using (B5), (B8) and the Mean Value Theorem  we can write
$$\big| b_{ij}(x_2,\zeta)\,-\,b_{ij}(x_2, t)\big|\,\leq\, \big| \p _r b_{ij}\big|_{C^1}\, \big|\,  \zeta \,-\, t\,\big|  \, \leq\, \,\kappa\,\epsilon\,+\,\kappa \, e^{^{1+\xi}}\,
$$  
\noindent And this together with (B5) imply 
$$\big| c  \big|_{C^0}\,\,\, \leq\, \,\,\kappa\,\epsilon\,+\,
 \kappa\, e^{^{1+\xi}}\, +\,   \big| b_{ij}(x,t)-e^t\delta_{ij}  \big|_{C^0} \,\,\,\leq\,\,\,\,\kappa\,\epsilon\,+\, \kappa\, e^{^{1+\xi}}\, +\,\epsilon=\,\,(1+\kappa)\,\epsilon\,+\, e^{^{1+\xi}}\, \kappa$$

\noindent {\bf The $C^1$-norm of $c$.} We have three types of first derivatives. First, from (B4) we have:
$$\p_t c\,=\,(\p_r d)\,(\p_t \zeta)\, +\, \big(   (\p_t \zeta)\,-\, 1    \big)\, e^\zeta\,I\, +\, (e^\zeta\,-\, e^t)\,I
$$
\noindent This last equation together with  (B5), (B6), (B7) imply
$$|\,\p_t c\,|\, \leq\, \epsilon \, (1+\kappa)\,+\,  \kappa\, e^{^{1+\xi}}\, +\, \kappa\, e^{^{1+\xi}}
= (1+\kappa)\,\epsilon \,+\,  2\, e^{^{1+\xi}}\, \kappa$$
\noindent where we are using the Mean Value Theorem, (B5)
and (B6) to estimate $e^\zeta-e^t$. Analogously
$$|\,\p_i c\,|\,=\,|(\,\p_r d)\,(\p_i \zeta)\, +\,    (\p_i \zeta)\, e^\zeta\, |\,\leq\, \epsilon\, \kappa\,+\, \kappa\,e^{^{1+\xi}}$$
\noindent and
$$|\,\barp_i c\,|\,=\,|\,\barp_i d |\,<\, \epsilon$$

\noindent {\bf The $C^2$-norm of $c$.} We have six types of first derivatives. As above using (B4),  (B5), (B6) and (B7)
we can obtain estimates for them. Here are the first three that do not involve the variable $x_2$:
{\small $$ \begin{array}{lll}
|\, \p^2_t c  \,|\,& =&\, \bigg|\,  (\p_r^2d)(\p_t\zeta)^2\,+\, (\p_rd)(\p_t^2\zeta)\, +\bigg[ (\p_t^2\zeta)\,
+\,\big((\p_t\zeta)^2\,-\,1\big)  \bigg]\,e^\zeta\,
 +\, (e^\zeta\,-\,e^t  )\, \bigg|\\ \\
& \leq&
\epsilon\,(1+\kappa)^2\, +\, \epsilon\,\kappa\,+\big[  \kappa+ 
\kappa(\kappa +2) \big]\,e^{^{1+\xi}}\,+\,\kappa\,e^{^{1+\xi}}\\ \\
&=& \big(1+3\kappa+\kappa^2\big)\,\epsilon\,+\, e^{^{1+\xi}}\big(4+\kappa \big)\,\kappa
\\ \\ \\
|\,\p_i\p_t c  \,|\,& =&\, \bigg|\,  (\p_i\p_rd)(\p_i\zeta)(\p_t\zeta)\,+\, (\p_rd)(\p_i\p_t\zeta)\, +
\bigg[ (\p_i\p_t\zeta)\,+\,(\p_i\zeta)(\p_t\zeta)  \bigg]\,e^\zeta  \, \bigg|\\ \\
& \leq\,&
\epsilon\,(1+\kappa)\,\kappa\, +\, \epsilon\,\kappa\,+\big[  \kappa+(1+\kappa)\kappa \big]\,e^{^{1+\xi}}\\ \\
&=&\big(2\kappa+\kappa^2\big)\,\epsilon\,+\,  e^{^{1+\xi}}\,
\big(2+\kappa\big)\,\kappa

\\ \\ \\
|\,\p_j\p_i c  \,|\,& =&\, \bigg|\,  (\p_j\p_id)(\p_j\zeta)(\p_i\zeta)\,+\, (\p_id)(\p_j\p_i\zeta)\, +
\bigg[ (\p_j\p_i\zeta)\,+\,(\p_j\zeta)(\p_i\zeta)  \bigg]\,e^\zeta  \, \bigg|\\ \\
& \leq\,&
\epsilon\,\kappa^2\, +\, \epsilon\,\kappa\,+\big[  \kappa+\kappa^2 \big]\,e^{^{1+\xi}}\\ \\
&=&
\big(\kappa+\kappa^2\big)\, \epsilon\,+\,e^{^{1+\xi}}\,\big(1+  \kappa\big)\,\kappa
\end{array}$$}\vspace{.1in}

\noindent And the ones involving the $x_2$ variable:
{\small$$
\begin{array}{lllll}
|\, \barp_j\barp_i c\,|\,& =&\,| \, \barp_j\barp_i d\,|\,\,\,\,<\,\,\,\,\epsilon&&\\ \\
|\, \p_j\barp_i c\,|\,& =&\,| \, (\p_r\barp_i d)\,(\p_j\zeta)\,|\,\,<\,\,\epsilon\,\kappa\\ \\
|\, \p_t\barp_i c\,|\,& =&\,| \,( \p_r\barp_i d)\,(\p_t\zeta)\,|\,\,<\,\,\epsilon\, (1+\kappa)\end{array}$$}
\vspace{.1in}

\noindent Note that all the estimates of the derivatives that we have 
obtained above are less or equal 
$$\begin{array}{lll}
\big(1+3\kappa+\kappa^2\big)\, \epsilon\,+\,e^{^{1+\xi}}\,\big(4+  2\kappa\big)\,\kappa&=&
\big(1+3\kappa+\kappa^2\big)\, \epsilon\,+\,e^{^{1+\xi}}\,\big(4+  2\kappa\big)\,L\,e^{-r\0{0}}\\\\
&\leq&\bigg[ 2(2+3\xi+\xi^2)e^{^{1+\xi}}L
\bigg]\big(\epsilon+e^{-r\0{0}}\big)\end{array}$$
\noindent And, since $\kappa\leq \xi$, we get
$|c|_{C^2}<2(2+3\xi+\xi^2)e^{^{1+\xi}}L(
\epsilon+e^{-r\0{0}})$. This concludes our calculations
and the proof of Claim 5.8.

Pedro Ontaneda

SUNY, Binghamton, N.Y., 13902, U.S.A.

\end{document}